\documentclass [reqno,11pt]{article}

\usepackage{amscd}
\usepackage[intlimits]{amsmath}
\usepackage{amssymb}

\newtheorem{teo}{Theorem}
\newtheorem{prop}[teo]{Proposition}
\newtheorem{lemma}[teo]{Lemma}
\newtheorem{example}[teo]{Example}
\newtheorem{remark}[teo]{Remark} 
\newtheorem{cor}[teo]{Corollary}
\newtheorem{definiz}[teo]{Definition}

\newtheorem{ack}{Acknowledgments}        

\newcommand{\fine}{\hfill$\Box$\\}
\newcommand{\dimo}[1][]         {\noindent\textbf{Proof#1}. }

\newcommand{\eps}               { \varepsilon}
\renewcommand{\phi}             {\varphi}
\newcommand{\id}                {\operatorname{\mathrm{id}}}

\newcommand{\restr}[1]          {\phantom{}_{\text{\raisebox{.4ex}{$|$}}#1}}
\newcommand{\R}                 { \mathbb { R } }
\newcommand{\C}                 { \mathbb { C } }
\newcommand{\Q}                 { \mathbb { Q } }
\newcommand{\Zeta}              { \mathbb {Z}  }

\newcommand{\I}                 { \operatorname{\mathrm {i} }}
\newcommand{\desu} [2] []       {\operatorname{\dfrac {\partial #1} {\partial #2}}}  
\newcommand{\desudt}[1] []      {\operatorname{\dfrac {\mathrm {d} #1 }{\mathrm {dt}}}}

\newcommand{\di}                {\, d }

\renewcommand{\det }            { \operatorname {det} }

\newcommand{\est}               {\raisebox{.34ex}
                                {$\scriptstyle {\bigwedge}$}}
\newcommand{\Aut}               {\operatorname {Aut}}

\newcommand{\vol}               {\operatorname{\mathrm {vol}}}
       
\newcommand{\Keler}             {K\"{a}hler }

\newcommand{\KE}                {K\"{a}hler-Einstein }

\newcommand{\debar }            {\bar {\partial } }

\newcommand{\Ric}{\operatorname {Ric} }

\newcommand{\dvol}              {\operatorname { dvol } }

\newcommand{\cinf}              {C^{\infty}}

\newcommand{\enf}               {\emph}
\newcommand{\Lap}               {\operatorname{\Delta}}
\newcommand{\OO}                {\mathcal{O}}
\newcommand{\chern}             {\operatorname{\mathrm{c}}}

\newcommand{\dphi}{\dot{\phi}}

\newcommand{\om}{\omega}
\newcommand{\omK}{{\omega_{KE}}}

\newcommand{\ra}                {\rightarrow}
\newcommand{\lra}               {\longrightarrow}
 
\newcommand{\unsu} [1]          { \dfrac {1} {#1}  } 
\newcommand{\menuno}            {{-1}}  
\newcommand{\vacuo}             {\emptyset}
\newcommand{\perdef }           { : \, = }
\newcommand{\idd}{\I \partial \debar}

\newcommand{\supp}{\operatorname{supp}}

\newcommand{\lamma}{\lambda}

\newcommand{\PP}{\mathbb{P}}
\newcommand{\alf}{\alpha}

\newcommand{\Y} { { (Y,\Delta_Y)} } 
\newcommand{\MV}{\langle [\om]^n, [X]\rangle}

\newcommand{\barPG}{P^0_G(\X,\om)} 
\newcommand{\PG}{P_G(\X,\om)}

\newcommand{\intv}{\frac{1}{V}\int_X} 
\newcommand{\Fo}{F^0}

\newcommand{\sing}{{\operatorname{sing}}}
\newcommand{\reg}{{\operatorname{reg}}}
\newcommand{\red}{{\operatorname{red}}}

\newcommand{\cordin}[1][z]{#1_1, ..., #1_n}

\newcommand{\galf}[1][f]{\operatorname{Gal}(#1)}
\newcommand{\galp}{\galf[\pi]}

\newcommand{\ord}{\operatorname{ord}}
\newcommand{\orb}{{{orb}}}
\newcommand{\mult}{\operatorname{mult}}

\newcommand{\XD}{\ensuremath{{(X,\Delta)} } }
\newcommand{\XDi}[1][i]{\ensuremath{{(X_{#1},\Delta_{#1})} } }

\newcommand{\xd}{X,\Delta}
\newcommand{\diam}{\operatorname{diam}}
\newcommand{\mc}[1][t]{\ensuremath{(*)_{#1} } }
\newcommand{\X}                 { \xd }

\begin{document}

\title{\KE metrics on orbifolds  and Einstein metrics on spheres}

\author{Alessandro Ghigi  and J\'anos Koll\'ar}

\maketitle

\begin{abstract}
\noindent
A construction of K\"ahler-Einstein metrics using Galois coverings,
studied by Arezzo--Ghigi--Pirola, is generalized to orbifolds.  By
applying it to certain orbifold covers of $\C\PP^n$ which are trivial
set theoretically, one obtains new Einstein metrics on odd-dimensional
spheres.  The method also gives K\"ahler-Einstein metrics on degree 2
Del Pezzo surfaces with $A_1$ or $A_2$--singularities.
\end{abstract}

%\normalsize

%\tableofcontents

\section{Introduction}

The aim of this paper is to explain how the methods of
Arezzo, Ghigi, and Pirola
\cite{arezzo-ghigi-pirola} can be applied to
construct \KE metrics on compact complex orbifolds
with positive first Chern class,
and then use the approach of Boyer, Galicki, and  Koll{\'a}r
\cite{boyer-galicki-kollar-Annals}
to obtain new Einstein metrics on odd dimensional spheres.

The somewhat unusual aspect  is that we work with
orbifolds $\mathcal X$ that admit a 
map $\pi:\mathcal X\to \PP^n$ which is the {\em identity} map
set theoretically. Nonetheless, in the orbifold category
$\pi$ is a nontrivial Galois cover, although with
trivial Galois group.
\medskip

The existence of \KE metrics on compact complex manifolds with
positive first Chern class is still a difficult problem.
For surfaces
and toric manifolds a complete solution is known, due respectively to
Tian \cite{tian-DP} and Wang-Zhu \cite{wang-zhu}. Apart from
these cases, there are two large classes of examples.  The simplest are
homogeneous spaces, for instance $\PP^n$, quadrics, Grassmannians.  In
all these cases, the first Chern class is large, meaning for instance,
that it is a large multiple of a generator of $H_2(X,\Zeta)$.  The
opposite case, when the first Chern class is a small multiple of a
generator of $H_2(X,\Zeta)$ is also understood in many instances; see
\cite{bourg-fano} for a good overview.

A blending of these two approaches was developed in 
Arezzo, Ghigi, and Pirola \cite{arezzo-ghigi-pirola} to yield \KE metrics on
certain manifolds $X$ which can be realized as Galois covers
of another manifold $Y$ with a \KE metric. Since the method relies on
finite group actions,  it is most successfull 
when symmetries form a natural part of the complex structure,
for instance for double covers of $\PP^n$.

A construction of Einstein metrics on
odd dimensional spheres was studied in
Boyer, Galicki, and  Koll{\'a}r 
\cite{boyer-galicki-kollar-Annals}. The idea is that
the quotient of an odd dimensional sphere by a circle action is
frequently a complex orbifold, and a result of
Kobayashi \cite{kob} allows one to lift a \KE orbifold metric from the quotient
to an Einstein metric  on the sphere.

A frequently occurring case, studied by Orlik and Wagreich 
\cite{or-wa} and Boyer, Galicki, and  Koll{\'a}r
\cite{boyer-galicki-kollar-Annals},
appears  when the quotient $S^{2n+1}/S^1$
is $\PP^n$  as a manifold, and  the orbifold structure
is given by 
a $\Q$-divisor 
$$
\Delta=\sum_{i=0}^{n+1}\bigl(1-\tfrac1{m_i}\bigr)D_i,
$$
where
$$
      D_i=\{z_i=0\} \quad\text { for }\quad i=0,...,n,\quad
      D_{n+1}=\{ z_0+\cdots +z_n=0\},
$$
and the $m_0,\dots,m_{n+1}$ are 
pairwise relatively prime ramification indices.
(See Section \ref{orb.section} for precise definitions.)
The orbifold first Chern class is
$$
c_1(\PP^n,\Delta)=(n+1)-\sum_{i=0}^{n+1} \bigl(1-\tfrac1{m_i}\bigr)=
\sum_{i=0}^{n+1}\tfrac1{m_i}-1,
$$
where we have identified $H^2(\PP^n,\Q)$ with $\Q$.
Thus $c_1(\PP^n,\Delta)$ is positive iff 
\begin{equation}
\sum_{i=0}^{n+1}\tfrac1{m_i}-1>0.
\label{c_1>0}
\end{equation}
The  existence result \cite[Theorem 34]{boyer-galicki-kollar-Annals}
shows that $(\PP^n,\Delta)$ has an orbifold  \KE metric if
in addition the following inequality is also satisfied
\begin{equation}
\sum_{i=0}^{n+1}\tfrac1{m_i}-1<\tfrac{n+1}{n}
\min_i \{     \tfrac{1}{m_i} \}.
\label{bgk.old}
\end{equation}
This paper started with the observation that
one can apply the method of
\cite{arezzo-ghigi-pirola}
to the {\em identity} map $(\PP^n,\Delta)\to \PP^n$
which is a Galois cover (with trivial Galois group).
On the other hand, over the affine chart
$\PP^n\setminus\{D_i\cup D_j\}$
the same map can be viewed as  having cyclic Galois group
of order $\prod_{k\neq i,j} m_k$.
This approach  improves the bound of 
\cite{boyer-galicki-kollar-Annals}
by a factor of $n$, and we obtain
\begin{teo} Let $D_0,\dots,D_{n+1}\subset \PP^n$ be hyperplanes in
  general position and $m_0,\dots,m_{n+1}$ pairwise relatively prime
  natural numbers.  Assume that
  \begin{equation}
    0<\sum_{i=0}^{n+1}\tfrac1{m_i}-1<(n+1)
    \min_i \{     \tfrac{1}{m_i} \}.
    \label{bgk.new}
  \end{equation}
  Then there is an orbifold \KE metric on
  $(\PP^n,\sum_{i=0}^{n+1}(1-\tfrac1{m_i})D_i)$.
\end{teo}
Set $M=\prod_i m_i$ and $w_i=M/m_i$.   
 As shown in \cite{boyer-galicki-kollar-Annals} the
intersection of the unit sphere with the Brieskorn--Pham singularity 
$$
L(m_0,\dots,m_{n+1}):=
S^{2n+3}\cap \bigl(\sum_{i=0}^{n+1} z_i^{m_i}=0\bigr)\subset \C^{n+2}
$$
is homeomorphic to  $S^{2n+1}$ and a \KE metric on the
corresponding projective orbifold 
$$
(X,\Delta_X):=\left(\bigl(\sum_{i=0}^{n+1} z_i^{m_i}=0\bigr),
\sum_{i=0}^{n+1}(1-\tfrac1{m_i})[z_i=0]\right)\subset
\PP(w_0, \dots, w_{n+1})
$$
 lifts to a
positive Ricci curvature Einstein metric on 
$L(m_0,\dots,m_{n+1})$.
The weighted projective space $\PP(w_0, \dots, w_{n+1})$
is not well formed  and it is isomorphic to
the ordinary projective space $\PP^{n+1}$ by the map
$$
(z_0,\dots,z_{n+1})\mapsto (x_0=z_0^{m_0},\dots, x_{n+1}=z_{n+1}^{m_{n+1}}).
$$
Under this isomorphism we get that
$$
(X,\Delta_X)\cong \left(\bigl(\sum_{i=0}^{n+1} x_i=0\bigr),
\sum_{i=0}^{n+1}(1-\tfrac1{m_i})[x_i=0]\right)\subset
\PP^{n+1}.
$$
By eliminating the variable $x_{n+1}$ we get that
$$
(X,\Delta_X)\cong (\PP^n,\Delta).
$$
The
isometry  class of the metric on the sphere determines the complex
orbifold $(\PP^n, \sum_{i=0}^{n+1}(1-\tfrac1{m_i})D_i)$,
except possibly when  $(\PP^n, \sum_{i=0}^{n+1}(1-\tfrac1{m_i})D_i)$
has a holomorphic contact structure. The latter can happen only when
$n$ is odd; see \cite[Lem.17]{boyer-galicki-kollar-Annals}
for another necessary condition.
(Note that
$n+2$ hyperplanes in general position do not have moduli, so the
numbers $m_0,\dots,m_{n+1}$ alone determine the complex orbifold.)

Even with the improved bounds, the equations
 (\ref{bgk.new}) are not easy to satisfy.
Still,  as in  Example \ref{new.exmps}, we get
 12 new Einstein metrics on $S^5$
corresponding to the ramification indices  
$$
m_0=2, m_1=3, m_2=5, 
m_3\in\{17,19,23,29,31,37,41,43,47,49,53,59\},
$$ 
$\geq 10^3$ new Einstein metrics on $S^7$,
$\geq 10^6$ new Einstein metrics on $S^9$ ...

The above construction can be varied in many ways.
For instance, one can take more than $n+2$ hyperplanes and quadrics.
In all of these cases one gets an improvement by a factor roughly
$n$ compared to the bounds in
\cite{boyer-galicki-kollar-Annals}, but this gives many new
cases only for $n$ large.
(As shown by Orlik and Wagreich \cite{or-wa}, taking higher degree 
hypersurfaces for the $D_i$ yields Einstein metrics
on various rational homology spheres.)

\medskip

As another application, we consider singular degree 2 Del Pezzo
surfaces.  These are all double covers of $\PP^2$ ramified along a
quartic curve.  In the smooth case the existence of \KE metrics was
proved by Tian \cite{tian-DP}.  For singular surfaces we get the following.

\begin{teo}\label{teo-dp}
  Let $S$ be a degree 2 Del Pezzo surface with only $A_1$ or $A_2$
  singularities. Then $S$ has an orbifold \KE metric.
\end{teo}

Anyone well versed in orbifolds, stacks and 
in the theory of 
Monge--Amp\`ere equations should have no problem
developing the theory of 
\cite{arezzo-ghigi-pirola}
in the orbifold setting.
Nonetheless, since the theory of orbifolds has too many
``well known'' but never proved theorems and not quite
correct definitions and proofs, we felt that it makes sense
to write down the arguments in some detail.

\section{Analytic coverings}

Let $X$ and $Y$ be reduced complex spaces.  A map $\pi : X \ra Y $ is
called \enf{finite} if it is proper and has finite fibres. Since $X$
is locally compact a finite to one map is proper if and only if it is
closed. Therefore a map is finite if and only if it is closed and has
finite fibres. (By contrast note that
$\pi:\C\setminus\{-1\}\to \{y^2=x^3+x^2 \}\subset \C^2$
given by $t\mapsto (t^2-1, t^3-t)$ is a closed map of algebraic
varieties with finite fibers but $\pi$ is not proper.)

The fundamental theorem on finite maps (see \cite[p.
179]{grauert-remmert-cas}) states that when $X$ and $Y$ are
irreducible any finite surjective map $\pi : X \ra Y$ is an
\enf{analytic covering}. This means that there is a thin subset
$T\subset Y$ such that
\begin{itemize}
\item [a)] $\pi^\menuno(T)$ is thin in $X$, and
\item [b)] the restriction $\pi^\menuno(Y\setminus T) \ra Y \setminus
  T$ is locally biholomorphic (\'etale).
\end{itemize}
Put $Y_0=Y\setminus T$ and $X_0=\pi^\menuno(Y_0)$. Then $\pi : X_0 \ra
Y_0$ is a topological covering. We call it a \enf{regular subcover} of
$\pi$.

We  assume that our spaces are irreducible so that
``analytic covering'' and ``finite holomorphic surjection'' can be
regarded as synonyms.

Another important fact is that an analytic covering $\pi : X \ra Y$
with $X$ and $Y$ normal is an open map (see \cite [p.
135]{grauert-remmert-cas}).

Let now $\pi : X \ra Y$ be an analytic covering among connected
\enf{normal} complex spaces. Put $Y'=\{y \in Y_\reg: \pi^\menuno(y)
\subset X_\reg\}$ and $X ' = \pi^\menuno(Y')$. Then $X'$ and $Y'$ are
open sets with complements of codimension at least 2. Now $\pi : X'
\ra Y'$ is a finite surjective map between complex manifolds.  Pick
local coordinates $z_1, ..., z_n$ on a neighbourhood $U$ of a point in
$X'$ and let $w_1, ..., w_n$ be coordinates around its image in $Y'$.
Let $w_i=\pi_i(z)$ be the local expression of $\pi$. The divisors
locally defined by the equation
\begin{equation*}
  \det \biggl ( \desu[\pi_i]{z_j} \biggr) =0
\end{equation*}
glue together yielding a well-defined divisor on $X'$. Since the
complement of $X'$ has codimension at least 2, the Remmert-Stein
extension theorem (see e.g. \cite[p. 181]{grauert-remmert-cas})
ensures that the topological closure of this divisor is a divisor in
$X$, called the \enf{ramification divisor} of $\pi$, and denoted by
$R=R(\pi)$.  It satisfies the Hurwitz formula $K_{Y'} = \pi^* K_{X'} +
R$.  Write $R=\sum_j r_j R_j$ with $R_j$ distinct prime divisors on
$X'$.  The reduced divisor $R_\red=\sum_j R_j$ is called the
\enf{ramification locus}.  By the
implicit function theorem $R_\red\cap X'$ is the set of points $x\in
X'$ such that $\pi$ is not \'etale at $x$, that is the set of critical
points of $\pi$. Since $\pi$ is finite, the image $\pi(R_\red)$ is a
divisor on $Y$, called the \enf{branch divisor} of $\pi$.

Consider now the sets $X''=X'\setminus \bigl ((R_\red)_\sing
\cup\pi^\menuno(B_\sing)\bigr )$ and $Y'= \pi(X'')$. Both are open and
have complements of codimension at least 2 in $X$ and $Y$
respectively. We use this notation often in the sequel. When we want
to stress the dependence on $\pi$, we write $X''(\pi)$ and $Y''(\pi)$.
If $x\in X''$ either $x\notin R_\red$ or $x$ belongs to one and only
one component $R_j$. In the first case we say that $\pi$ is
\enf{unramified} at $x$, in the latter case we say that the
\enf{ramification order of $\pi$ at $x$} is $r_j +1$. The ramification
order of $\pi$ at $x$ will be denoted by $\ord_\pi(x)$. When $\pi$ is
unramified at $x$, we put $\ord_\pi(x)=1$.  If $D\subset X$ is an
irreducible divisor, then there is an open dense subset $D''\subset D$
such that $\ord_\pi(x)$ does not depend on $x\in D''$. This common
value is denoted by $\ord_\pi(D)$ and it is called the \enf{
  ramification order of $\pi$ along $D$}.

We use some basic properties of analytic coverings and maps between
them (see, for instance, \cite[Lemma 16.1]{barth-peters-vandeven}).

\begin{lemma}\label{BPV}
  Let $x\in X''$. If $\pi$ is unramified at $x$, then $\pi$ is a local
  biholomorphism at $x$. If it has ramification order $m>1$, let $R_j$
  be the component of $R_\red$ passing through $x$. Then there are
  local coordinates $\cordin[z]$ on $X''$ and $\cordin[w]$ on $Y''$
  centred at $x$ and $y=\pi(x)$ respectively, such that locally $R_j=
  \{z_1=0\} $, $B=\{w_1=0\}$ and $\pi(z_1, ..., z_n) = (z_1^{m} , z_2,
  ..., z_n)$.
\end{lemma}
Since the complement of $X''$ has codimension 2, $R_\red$ is the
closure of $R_\red\cap X''$, that is the closure of the set of points
where $\pi$ has ramification order $>1$.

The next lemma considers the problem of lifting in the simplest case.
Denote by $D(r)$ the disc of radius $r$ centred at the origin, by
$D^*(r)$ the complement of $\{0\}$ in $D(r)$, and by $P(r_1, ..., r_n)
$ the polydisc centred at the origin with polyradius $(r_1, ...,
r_n)$.
\begin{lemma}
  Let $ P_1= P(r_1, ..., r_n), P_2=P(\rho_1, ...,\rho_n)$, $Q_1=
  P(r^{m_1}_1, r_2, ..., r_n)$, $Q_2=P(\rho^{m_2}_1, \rho_2,
  ...,\rho_n )$. Set $P_1^*=D^*(r_1)\times P(r_2, ..., r_n)$ and
  similarly for $P_2^*, Q_1^*, Q_2^* $.  Let $\pi_i:P_i \ra Q_i$ be
  the maps $ \pi_1(z_1,..,z_n) = (z_1^{m_1},z_2,...,z_n) $, $
  \pi_2(z_1,..,z_n) = (z_1^{m_2},z_2,...,z_n) .  $ Let $f: Q_1 \ra
  Q_2$ be a holomorphic map such that $f(Q_1^*) \subset Q_2^*$.  If
  $m_2|m_1$ there are exactly $m_2$ liftings of $f$ (that is maps
  $\tilde {f}: P_1\ra P_2$ such that $\pi_2 \tilde{f}=f\pi_1$).  Any
  local lifting of $f$ defined in a neighbourhood of some point $x\in
  P_1$ extends to one of these liftings defined on $P_1$.
\end{lemma}

\begin{lemma}\label{monodromy.lem}
  Let $\pi_1: X_1 \ra Y$ and $\pi_2 :X_2 \ra Y$ be analytic
  coverings.  For $U\subset
  X_1$ set
  \begin{equation*}
    \mathfrak{F}(U)=\{\text{holomorphic maps } s: U \ra X_2 \text{ such
      that } \pi_1 = \pi_2\circ s\}.
  \end{equation*}
  Then $\mathfrak{F}$ is a Hausdorff sheaf (of sets) over $X_1$.
  Assume that for any $x_1\in X''_1, x_2\in X''_2$ with $\pi_1(x_1) =
  \pi_2(x_2)$
$$
\ord_{\pi_2}(x_2) |
\ord_{\pi_1}(x_1) . 
$$
Then the restriction of $\mathfrak{F}$ to $X''_1\cap \pi_1^\menuno
Y''_2(\pi_2)$ is a finite topological covering.  In particular, if
$X''_1$ is simply connected, then for every $x_1\in X''_1\cap
\pi_1^\menuno Y''_2(\pi_2)$ and $x_2\in X''_2$ such that
$\pi_1(x_1)=\pi_2(x_2)$ there is an analytic map $f:X''_1\to X_2$ such
that $f(x_1)=x_2$ and $\pi_1=\pi_2\circ f$.
\end{lemma}

In fact, the above $f$ extends to $X_1$ by the following
immediate consequence of the Riemann Extension Theorem (see e.g.
\cite[p.144]{grauert-remmert-cas})

\begin{lemma}\label{riem.ext.lem}
 Let $\pi_1: X_1 \ra Y$ and $\pi_2 :X_2 \ra Y$ be analytic
  coverings, $X_1$ normal  and $T\subset X_1$ a thin set.
Let $f^0:X_1\setminus T\to X_2$ be an analytic map
such that $\pi_1=\pi_2\circ f^0$. Then $f^0$
extends to $f:X_1\to X_2$ such that 
$\pi_1=\pi_2\circ f$.
\end{lemma}

\section{The Galois group of  coverings}

Let $f:X\to Y$ be an analytic covering of normal complex spaces.  Put
$\galp=\{ f \in \Aut(X) : \pi \circ f = \pi\}$.  $\galp$ is a finite
subgroup of $\Aut(X)$.  In fact fix $x\in X''$, $y=\pi(x)$, and let
$V$ be a neighbourhood of $y$ in $Y$ such that $\pi^\menuno(V) =
\bigcup_{i=1 } ^k U_i$ with $\pi: U_i \ra V$ a biholomorphism and
$x\in U_1$. Then the stabiliser $\galp_x$ is a subgroup of finite
index in $\galp$.  Moreover any $f\in \galp_x$ maps $U_1$ to itself.
Since $\pi\restr{U_1}$ is injective, the restriction of $f$ to $U_1$
is the identity. By the connectedness of $X$, $f=\id_X$, so
$\galp_x=\{1\}$ and $\galp$ is finite.

Since $\pi$ is $\galp$-invariant, the $\galp$-orbit of $x\in X$ is
contained in $\pi^\menuno\bigl (\pi(x) \bigr )$.  We say that an
analytic covering $\pi : X \ra Y$ is \enf{Galois} if the converse
holds, that is two points of $X$ lie on the same fibre of $\pi$ only
if they belong to the same $\galp$-orbit.

The branching divisor of a Galois cover can be described also in the
following way.  Given a prime divisor $D$ in $X$, set
$\Gamma(D)=\{\gamma\in \Gamma: D \subset
\operatorname{Fix}(\gamma)\}$.  For each prime divisor $D$ the image
$\pi(D)$ is a prime divisor in $Y$.  The prime divisors for which
$\Gamma(D)\neq 0$ are exactly the $R_j$. Set $B_j=\pi(R_j)$. In
general different $R_j$'s can have the same image. Assume that
$\{B_i\}_{i\in I}$ is the set of all images of the $R_j$'s (that is
$B_i\neq B_k$ if $i\neq k$). Then
\begin{equation}
  B(\pi) = \sum_{i\in I}\biggl ( 1 - \frac{1}{|\Gamma(R_i)|}\biggr) B_i.
\end{equation}

\begin{lemma}\label{regular-Galois}
  Let $X$ and $Y$ be normal complex spaces, $\pi : X \ra Y$ an
  analytic covering and $Y_0\subset Y$ an open subset with thin
  complement.  Put $X_0=\pi^\menuno (Y_0)$ and $\pi_0=\pi\restr{X_0} :
  X_0 \ra Y_0$. Then the elements of $\galf[\pi_0]$ extend to elements
  of $\galp$, and if $\pi_0$ is Galois, then $\pi$ is Galois too.
\end{lemma}
\dimo The first part follows from Lemma (\ref{riem.ext.lem}).  For the
second part, let $x, x'\in X$ be such that $\pi(x)=\pi(x')=y$. If
$y\in Y_0$ there is some $g\in \galf[\pi_0]$ such that $g.x=x'$. Since
we have just proved that $\galf[\pi_0]=\galp $ the Galois condition is
satisfied for these points. If instead $y\in Y \setminus Y_0$, choose
neighbourhoods $U_i$ and $V$ as above. Assume $x=x_1 \in U_1$ and
$x'=x_2\in U_2$.  Let $\{z_n\}$ be a sequence of points in $X_0\cap
U_1$ converging to $x$.  Then $y_n=\pi(z_n)$ converge to $y$. Since
$\pi$ is open, $\pi(U_2)=V$. Therefore there are points $z'_n \in
U_2\cap X_0$ such that $\pi(z'_n)=y_n$. By the Galois condition on
$X_0$, there are $g_n \in \galp$ such that $z'_n = g_n . z_n$. As
$\galp$ is finite, we can extract a subsequence with $g_n \equiv g$.
Since $\lim z'_n = x_2 $ as $\pi^\menuno(y)\cap U_2 =\{x_2\}$, we get
$x_2 = g . x_1$.  \fine

If $\pi : X \ra Y$ is a Galois covering, then $\galp$ acts freely on
any regular subcover $X_0$. Therefore if $x, x'\in X_0$ 
and  $\pi(x)=\pi(x')$, then  there is a
unique $g\in \galp$ such that $g.x=x'$.  In particular the cardinality
of $\galp$ equals that of the generic fibre. This condition is also
sufficient: $\pi$ is Galois iff $|\galp|$ equals the cardinality of
the general fibre iff $\galp$ is transitive on the general fibre.

For later reference we state the following simple
lemma.

\begin{lemma}\label{Galois-triple}
  Let $X, Y$ and $ Z$ be irreducible complex spaces, and $f: X\ra Z$,
  $g: Y\ra Z$, $h: X \ra Y$ analytic coverings such that $gh=f$. If $f$
  is Galois, then $h$ is Galois too.
\end{lemma}
\dimo Thanks to Lemma \ref{regular-Galois} it is enough to consider
the unramified case. Fix $x\in X$ and put $y= h(x), z= f(x)=g(y)$. We
need to show that $h_*\pi_1(X,x) $ is a normal subgroup of
$\pi_1(Y,y)$. Since $g_* : \pi_1(Y,y) \ra \pi_1(Z,z)$ is injective it is
enough to check that $g_* h_*\pi_1(X,x) $ is a normal subgroup of
$g_*\pi_1(Y,y)$. But $f$ being Galois $f_*\pi_1(X,x) = g_*
h_*\pi_1(X,x)$ is normal in $\pi_1(Z,z)$, hence a fortiori in $g_*
\pi_1(Y,y)$.  \fine

For a general analytic covering $\pi : X \ra Y$ it is not possible to
assign multiplicity to the branching divisor in any reasonable way. In
fact, different points in the preimage of a point $y\in B$ have
different branching orders. A typical example is 
$X=\{z^3-3yz+2x=0\}\subset \C^3$ projecting on $\C^2_{x,y}$.
Even shrinking the domain around the origin, 
one cannot separate the branches with
different orders.

On the other hand, when the covering is Galois, for any $y\in Y''$ all
points in $\pi^\menuno(y)$ have the same branching order. Therefore we
can assign multiplicities to the branch divisor according to the
following rule. Let $y\in Y''\cap B$ and let $x$ be any point in $
\pi^\menuno(y)$. Then we define the multiplicity of $B$ in $y$ to be $
1-1/\ord_{\pi}(x)$. We still denote by $B$ the $\Q$-divisor given by the
branching locus provided with these multiplicities.
Note that with this convention $R=\pi^*B$, that is, the ramification
divisor is the pull back of the branch divisor.

\section{Orbifolds as pairs}
\label{orb.section}

As in  \cite{boyer-galicki-kollar-Annals}, we look at orbifolds as a
particular type of log pairs. $(X, \Delta)$ is a log pair if $X$ is a
normal algebraic variety (or a normal complex space) and 
$\Delta=\sum_i d_i D_i$ is
an effective $\Q$-divisor where the $D_i$ are distinct, irreducible divisors
and $d_i\in\Q$. The number $d_i$ is called
the  multiplicity of $\Delta$  along $D_i$, it is denoted by
$\mult_{D_i}\Delta$. We set $\mult_{D}\Delta=0$ for every 
other irreducible divisor $D\neq D_i\ \forall i$.

 Let $X''(\Delta)$ (or simply $X''$) be
the complement of $X_\sing \cup \Delta_\sing$. For $x\in X''$ the
multiplicity of $\Delta$ at $x$ is a well defined rational number.  
For orbifolds, we need to 
 consider only pairs $(X,\Delta)$ such that $\Delta$ has the form
\begin{equation*}
  \Delta=\sum_i \bigl ( 1 - \tfrac{1}{m_i} \bigr ) D_i,
\end{equation*}
where the  $D_i$ are prime divisors and $m_i \in \mathbb{N}$.  If $(X,\Delta)$
is such a pair then for any divisor $D\subset X$ we put
\begin{equation*}
  \ord_{\Delta}(D) = \frac {1}{1-\mult_{D}\Delta}.
\end{equation*}
The
assumption on the multiplicities of $\Delta$ amounts to saying that
the order is always a nonnegative integer.  
\begin{definiz}
  An \enf{orbifold chart} on $X$ compatible with $\Delta$ is a Galois
  covering $\phi : U \ra \phi (U) \subset X$ such that 
\begin{enumerate}
\item $U$ is a
  domain in $\C^n$ and $\phi(U)$ is open in $X$; 
\item  the branch locus
  of $\phi$ is $\Delta_\red \cap \phi(U)$; 
\item  for any $x\in
  U''(\phi)$ such that $\phi(x)\in D_i$, $\ord_\phi(x)=m_i$.
 \end{enumerate}
 Conditions (2) and (3) are equivalent to 
  \begin{equation}
    B(\phi)= \Delta\cap \phi(U). 
  \end{equation}
\end{definiz}
\begin{definiz}
  An orbifold is a log pair $(X,\Delta)$ such that $X$ is covered by
  orbifold charts compatible with $\Delta$.
\end{definiz}
(For a slightly more general approach, see \cite[\S
14]{deligne-mostow}.)

Let $X$ be a normal complex space and $\pi:U\to X$ a Galois cover
where $U$ is a smooth. As discussed earlier, the branch divisor
$B(\pi)$ of $\pi$ is defined and we get a log pair $(X, B(\pi))$. If
$U$ is simply connected, (which we can always assume by shrinking $U$
suitably) then by Lemma \ref{monodromy.lem} the log pair $(X, B(\pi))$
determines $\pi:U\to X$ up to biholomorhisms. Thus we recover the
classical definition of orbifolds (as in \cite{baily-decomposition}
for example).

\begin{example}\label{local-normal-crossing}
  Let $X$ be a complex manifold and $D=\sum_{i\in I} D_i$ a divisor
  with \enf{local} normal crossing. By this we mean that for any point
  $x\in X$ there is a holomorphic coordinate system $(V,z_1, ...,
  z_n)$ such that $D\cap V = \{z\in V: z_1\cdots z_k=0\}$. If $D_i\cap
  V\neq \emptyset$ then $D_i\cap V$ is the union of some of the
  hypersurfaces $\{z_j=0\}$. ($D$ is said to be a divisor with
  \enf{global} normal crossing if, in addition,  each  $D_i$ is
 smooth.)  For any $i\in I$, fix an integer $m_i >1$ and put
  $\Delta=\sum_i ( 1- 1/m_i) D_i$. We claim that $(X,\Delta)$ is an
  orbifold.  Indeed, fix a coordinate system as above and put $m'_j =
  m_i$ if $\{z_j=0\} \subset D_i\cap V$. Set
  \begin{equation}
    \label{eq:snc-chart}
      \phi: U\ra V, \qquad \pi(x_1, ..., x_n) = (x_1^{m'_1}, ...,
      x_k^{m'_k}, x_{k+1}, ..., x_n).
  \end{equation}
  Then $(U,\phi)$ is orbifold chart on $X$ compatible with $\Delta$
  and so $\XD$ is an orbifold.
\end{example}

In the same way, the usual definition of orbifold map is equivalent to
the following one.
\begin{definiz}
  A finite holomorphic map $f: X \ra Y$ is an \enf{ orbifold map} $f:
  (X,\Delta_X) \ra (Y, \Delta_Y)$ if
  \begin{equation}
    \label{eq:orbimap}
    \ord_{\Delta_Y}(f(D)) \Big | \ord_{\Delta_X}(D) \cdot \ord_{f} D    
  \end{equation}
  for every divisor $D\subset X$.

An  \enf{orbifold automorphism} is an orbifold map that is invertible
  with inverse an orbifold map. The group of automorphisms of
  $(X,\Delta)$  is denoted by $\Aut(X,\Delta)$.
\end{definiz}

\begin{definiz}
\label{orbifold-Galois}
  An orbifold Galois covering $f : (X,\Delta_X) \ra (Y,\Delta_Y)$ is
  an orbifold map such that $f: X \ra Y$ is a Galois analytic cover
  and $\galf \subset \Aut(\xd_X)$.
\end{definiz}
By the \enf{degree} of an orbifold Galois cover we mean its
\emph{degree} as an analytic cover.

\begin{lemma}\label{Galois-liftings}
  Let $f : (X,\Delta_X) \ra (Y,\Delta_Y)$ be an orbifold map.  Then
  given $x\in X$ and $y=f(x)\in Y$ there are orbifold charts
  $(U,\phi)$ and $(V,\psi)$ around $x$ and $y$ respectively such that
  $f$ has a lifting $\tilde{f}: U \ra V$.  If, in addition, $f : X\ra
  Y$ is a Galois covering then $\tilde{f}: U \ra V$ is also a Galois
  covering.
\end{lemma}
\dimo Choose the chart $(U,\phi)$ such that $U$ is simply connected
and $f\bigl ( \phi (U)\bigr )\subset \psi(V)$.  If $D\subset U$ is any
divisor then
$$
\ord_{f\circ \phi}D=\ord_f{\phi}(D)\cdot \ord_{\phi}D
=\ord_f{\phi}(D)\cdot \ord_{\Delta_X}\phi(D).
$$
By the definition of orbifold maps,
$$
\ord_{\Delta_Y}(f\circ\phi)(D) \Big| \ord_f{\phi}(D)\cdot \ord_{\Delta_X}\phi(D),
$$
hence we conclude that $\ord_{\Delta_Y}(f\circ\phi)(D)$ divides $
\ord_{f\circ \phi}D$.  Thus the assumption of Lemma
\ref{monodromy.lem} is satisfied and so $f\circ \phi$ lifts to $\tilde
f:U\to V$.  Assume next that $f : X\ra Y$ is a Galois covering and
pick $u_1,u_2\in U$ such that $\tilde f(u_1)=\tilde f(u_2)$. Then
$f(\phi(u_1))=f(\phi(u_2))$ hence there is a Galois automorphism
$\sigma$ of $f$ such that $\phi(u_1)=\sigma(\phi(u_2))$.  Applying
Lemma \ref{monodromy.lem} to $\phi:U\to X$ and $\sigma\circ \phi:U\to
X$ we conclude that $\sigma$ lifts to a biholomorphism $\tilde \sigma$
of $U$ such that $\phi(u_1)=\phi(\tilde \sigma(u_2))$.  Since
$\phi:U\to X$ is Galois, there is a biholomorphism $\rho$ of $U$ such
that $u_1=\rho(\tilde \sigma(u_2))$.  This shows that in the
commutative diagram
\begin{equation}\label{eq:charts-lifting}
  \begin{CD}
    U @>\tilde{f}>> V \\
    @VV{\phi}V @VV{\psi}V\\
    \phi(U) @>f>> \psi(V).
  \end{CD}
\end{equation}
the composite $f\circ \phi$  is Galois. But
$f\phi= \psi \tilde{f}$ and by Lemma
\ref{Galois-triple} $\tilde{f}$ is a Galois cover.  \fine
\begin{example}\label{snc-galois}
  Let $(X,\Delta)$ be any orbifold, and let $(X,0)$ denote the
  orbifold structure on $X$ with trivial branching divisor.  It is a
  nontrivial result that $(X,0)$ is an orbifold, that is, $X$ has
  quotient singularities (see \cite{prill}).  (We use mainly the case
  when $X$ is smooth, and then the orbifold charts of $(X,0)$ are
  simply the manifold charts of $X$.)

  The identity map $\id_X : \XD \ra (X,0)$ is trivially an orbifold
  Galois covering. In fact it is both an orbifold map and a Galois
  analytic cover, and $\galf[\id_X]=\{id_X\}\subset \Aut\XD$.
\end{example}
If $f : (X,\Delta)\ra (Y,\Delta_Y)$ is an orbifold Galois covering the
\enf{orbifold ramification divisor} of $f$ is defined as
\begin{equation*}
  R^\orb(\Delta_X, \Delta_Y,f) = R(f) + \Delta_X - f^* \Delta_Y.  
\end{equation*}
With this definition the logarithmic ramification formula
\begin{equation*}
  K_X + \Delta_X = f^* ( K_Y + \Delta_Y) + R^\orb (\Delta_X, \Delta_Y, f) 
\end{equation*}
is automatically satisfied. To understand the geometric meaning of
$R^\orb$ it is useful to look at the open set
\begin{equation*}
  X''(\Delta_X, \Delta_Y, f) = X_\reg \cap f^\menuno\bigl(Y_\reg
  \setminus (\Delta_Y \cup B(f) )_\sing\bigr) \setminus (\Delta_X \cup
  R(f) )_\sing.
\end{equation*}
This means that $x\in X''= X''(\Delta_X, \Delta_Y, f)$ if (a) $X$ is
smooth at $x$, (b) $Y$ is smooth at $y=f(x)$, (c) $x$ belongs to at
most one component $D$ of $\Delta_X + R(f)$ and in this case $x$ is a
smooth point of $D$, (d) $y$ belongs to at most one component $D'$ of
$\Delta_Y + B(f)$ and in this case it is a smooth point of $D'$.  As
usual the complement of this set has codimension 2.  Let $D$ be any
smooth divisor passing through $x$ and $D'$ a smooth component passing
through $y$.  Assume first that $f$ is unbranched at $x$ and that
locally $\Delta_X= (1- 1/p) D$ and $\Delta_Y=(1-1/q) D'$. Then there
is a local diagram like \eqref{eq:charts-lifting}, with $p=\deg \phi$
and $q=\deg\psi$.  Put $k=\deg\tilde{f}$. Since $f$ is unbranched we
can assume that its restriction to $\phi(U)$ is a biholomorphism onto
$\psi(V)$. Therefore $p=qk$. If $p=1$, then $q=k=1$, and as expected
$\mult_x R^\orb=0$.  If $p>1$, then necessarily $D'=f(D)$ because of
\eqref{eq:orbimap} and $f^* D'= D$, since $f$ is \'etale. Therefore
$R^\orb = (1/q-1/p) D = (k-1)/p \cdot D$.  If instead $\ord_x(f)=m
>1$, then again $D'=f(D)$, $R(f)=(m-1)D$, $f^*D'=mD$, $pm=qk$ and
$R^\orb= (m/q - 1/p) D = (k-1)/p \cdot D$ once more.  Roughly the
orbifold ramification divisor is the ramification of the lifting
$\tilde{f}$ divided the degree of the local chart $\phi$.

Let $(X,\Delta)$ be an orbifold and $\Gamma \subset \Aut(X,\Delta)$ a
finite subgroup. We want to define a quotient orbifold $(Y,\Delta')$.
By Cartan's lemma \cite{cartan-pro-lefschetz-suo} $Y=X/\Gamma$ is a
normal analytic space and the canonical projection $\pi: X\ra Y$ is an
analytic covering.  The support of the branch divisor $\Delta'$ is
defined to be $\pi(\Delta)\cup B(\pi)$, while the multiplicities are
specified as follows.  Let $D$ be an irreducible component of
$\pi(\Delta)\cup B(\pi)$.  If $D$ is a component of $\pi(\Delta)$ and
not of $B(\pi)$, then we assign to $D$ the multiplicity $
\mult_x(\Delta)$, where $x$ is any point in $ X''(\Delta)$ such that
$\pi(x) \in D $ is a smooth point of $ \pi(\Delta) \cup B(\pi)$.  If
$D$ is a component of $B(\pi)$ and not of $\pi(\Delta)$ then we assign
to $D$ the same multiplicity it has as a component of $B(\pi)$, that
is $1-1/\ord_\pi(x)$ for any $x\in X''(\pi)$ such that $\pi(x) \in D $
is a smooth point of $ \pi(\Delta)\cup B(\pi)$.  Finally, if $D$ is a
common component of $\pi(\Delta)$ and $B(\pi)$ the we assign to it the
multiplicity
\begin{equation*}
  1 - \frac{1-\mult_x\Delta}{\ord_\pi(x)}
\end{equation*}
for any $x\in X''(\Delta) \cap X''(\pi)$ such that $\pi(x) \in D $ is
a smooth point of $ \pi(\Delta)\cup B(\pi)$.

\begin{prop}
  Let $(X,\Delta)$ be an orbifold, and $\Gamma\subset \Aut(\xd)$ a
  finite subgroup.  Let $Y=X/\Gamma$ be the quotient analytic space,
  and $\Delta'$ the $\Q$-divisor defined above. Then $(Y,\Delta')$ is
  an orbifold and the canonical projection
  \begin{equation}
    \pi : (X, \Delta_X) \lra (Y, \Delta')
  \end{equation}
  is an orbifold Galois covering.
\end{prop}
\dimo We need to show that $Y$ is covered by orbifold charts
compatible with $\Delta'$. Fix $y\in Y$, $x\in \pi^\menuno(y)$ and let
$\phi: U \ra \phi(U)$ be an orbifold chart with $x\in \phi(U)$.  If
the stabiliser $\Gamma_x$ is trivial we can assume that $\gamma
\phi(U) \cap \phi(U) =\vacuo$ for any $\gamma\neq e$.  Then $\pi:
\phi(U) \ra Y$ is a biholomorphism onto its image. Put $\psi=\pi\phi :
U \ra Y$. We claim that $\psi$ is an orbifold chart on $Y$ compatible
with $\Delta'$. In fact $\psi$ is Galois since $\pi$ is a
biholomorphism on $\phi(U)$, and $\pi^*B(\psi) = B(\phi)= \Delta\cap
\phi(U)$.  On the other hand $B(\pi)\cap \psi(U) =\vacuo$ since $\pi :
\phi(U) \ra \psi(U)$ is biholomorphic. Therefore on $\psi(U) $ the
divisor $\Delta'$ coincides with $B(\psi)$. This proves that $\psi: U
\ra Y$ is an orbifold chart.  If $\Gamma_y \neq \{e\}$ take a chart
$\phi: U \ra \phi(U)\subset X$ such that $\phi(U)$ be a
$\Gamma_x$-invariant neighbourhood of $x$. Lemma
\ref{Galois-liftings} ensures that also in this case $\psi=\pi \phi: U
\ra \psi(U) \cong \phi(U)/\Gamma_x$ is a Galois covering. It is easy
to verify that $B(\psi)= \Delta'$ on $\psi(U)$.  Finally that $\pi$ is
an orbifold Galois covering is clear: a lifting of $\pi : \phi(U) \ra
\psi(U)$ is given by the identity map $U\ra U$ which is trivially
Galois.  \fine

\section{Basic estimates for orbifold \KE metrics}

In this section we collect the orbifold versions of some fundamental
results due to  Aubin, Bando-Mabuchi and Tian,
that are needed in the existence criteria in the next section. Most of
the proofs are the same as in the case of a manifold and we just give
appropriate references.  For the basic definitions of differential
geometry on orbifolds see \cite{baily-decomposition},
\cite{baily-imbedding}, \cite{boyer-galicki-twistor} and
\cite{borzellino-indiana} . Some information on Sobolev spaces and
Laplace operators on orbifolds can be found e.g. in \cite{chiang-AMS}.

\begin{remark}
  Note  that if $X$ is a complex manifold and
  $\Delta$ is a non trivial branching divisor, then smoothness in the
  orbifold sense is rather different from ordinary smoothness.  For
  example, $f(z)=|z|$ is not smooth in the ordinary sense, but it
  belongs to $ \cinf(\C, \Delta)$, where $\Delta$ is the divisor
  concentrated at the origin with multiplicity $1/2$.  In fact the
  inclusions $\cinf(X) \subsetneq \cinf\XD$ and $\est^k(X) \subsetneq
  \est^k(\xd)$ are in general strict.
\end{remark}

\begin{definiz}
  A \enf{Fano orbifold} is a compact complex orbifold \XD such that
  $-(K_X + \Delta)$ is ample.
\end{definiz}
By the Baily-Kodaira imbedding theorem \cite {baily-imbedding} this is
equivalent to the fact that $\chern_1\XD$ contains an orbifold \Keler
metric.

The following is the orbifold analogue of Bonnet-Myers Theorem. It
follows, for example, from the Bishop volume comparison Theorem for
orbifolds, see \cite[Prop. 20, Cor. 21]{borzellino-indiana}.
\begin{teo}
  Let $X$ be an $m$\--dimensional orbifold and $g$ a Riemannian
  orbifold metric on $X$ with $\Ric(g)\geq \eps(m-1) g$ for some
  $\eps>0$. Then $\diam(X,g) \leq \pi/\sqrt{\eps}$.
\end{teo}

\begin{teo}
  [\protect{\cite[Theorem B]{nakagawa-isoperimetrico}}] Let $(X,g)$ be
  a Riemannian orbifold of dimension $m>2$ with $\Ric(g)\geq
  -(m-1)\eps^2 g$ for some $\eps\geq0$. Then there is a constant $C>0$
  depending only on $m$ and $\eps\cdot \diam(X,g)$ such that
  \begin{equation}
    \label{eq:gallot-li}
    ||\nabla u ||_{L^2} \geq
    C \frac{\vol(X,g)^{1/m}}{\diam(X,g)}
    ||u||_{L^{2m/(m-2)}}
  \end{equation}
  for any $u\in W^{1,2}(X) $ with $\int_X u \dvol_g =0$.
\end{teo}
Combining the last two theorems one gets the following uniform Sobolev
embedding.
\begin{cor}\label{uniform-Sobolev}
  Let $(X, \Delta)$ be an $n$\--dim\-ensional Fano orbifold.  For any
  $\eps>0$ there is a constant $C=C(\eps)>0$ such that for any metric
  $\om$ in the class $2\pi\chern_1(X,\Delta)$ with $\Ric(\om)\geq \eps
  \om$ and any $u\in W^{1,2}(X,\Delta)$
  \begin{equation}
    ||u||_{L^{2n/(n-1)}} \leq C 
    ||u||^2_{W^{1,2}}.
  \end{equation}
\end{cor}

If $(X,\Delta)$ is a \Keler orbifold, $\om\in \est^{1,1}(\xd)$ is a
closed smooth form and $\phi \in \cinf(\xd)$, put $\om_{\phi}=\om +
\idd \phi$. We write $\om_\phi >0$ to mean that it is a \Keler metric.
If $\om $ is such that
\begin{equation*}
  \langle [\om]^n , [X]\rangle = \int_X \om^n >0 
\end{equation*}
and $\phi\in\cinf\XD$ put
\begin{align}
  I_\om(\phi)&=\frac{1}{\MV}\int \phi (\om^n - \om^n_\phi)
  \label{def-di-I}
  \\ J_\om(\phi) & = \int_0^1 \frac{I_\om(s\phi)}{s}\di s
  \label{eq:def-di-J}
  \\
  \Fo_\om(\phi) & = J_\om(\phi) - \frac{1}{\MV}\int \phi \om^n.
\end{align}

\begin{lemma}
%  Let $\XD$, $\om$ and $\phi$ be as above. Then
  \begin{gather}
    \label{eq:J-stile-Tian-per-parti}
    J_\om(\phi) = \frac{1}{\MV} \sum_{k=0}^{n-1} \frac{k+1}{n+1}
    \int_M \I
    \partial\phi \wedge\debar \phi \wedge \om^k \wedge
    \om_\phi^{n-k-1}\\
    I_\om(\phi) -J_\om(\phi) =\unsu{\MV} \sum_{k=0}^{n-1}
    \frac{n-k}{n+1} \int_X \I \partial \phi \wedge \debar \phi \wedge
    \om^k \wedge \om_\phi^{n-k-1}.
    \label{eq:I-J}
  \end{gather}
  If $\om>0$ and $ \om_\phi>0$, then $I_\om(\phi)$, $J_\om(\phi)$ and
  $I_\om(\phi) - J_\om(\phi)$ are nonnegative and vanish only if
  $\phi$ is constant.  Moreover $J_\om \leq I_\om \leq (n+1)J_\om$.
\end{lemma}
For \eqref{eq:J-stile-Tian-per-parti} see \cite[Lemma
2.2]{tian-certain} or \cite[Lemma 2.1]{arezzo-ghigi-pirola}.  For
\eqref{eq:I-J} expand $\om^n - \om_\phi^n$.  The last statements
follow diagonalising simultaneously $\om$ and $\idd \phi$.  \fine

\begin{lemma}\label{lemma-triciclico}
 % Let $\XD$ and $\om$ be as above. 
  If $\lambda$ is a positive constant then
  \begin{equation}
    \label{eq:scaling-di-F0}
    F^0_{\lamma\om}(\lamma \phi) = \lamma F^0_\om( \phi). 
  \end{equation}
  Let $\om_0$ be a closed (1,1)-form such that $\langle [\om_0]^n,
  [X]\rangle >0$.  Given $\phi_{01}$, $ \phi_{12}\in \cinf\XD$ put
  $\om_1=\om_0 + \idd \phi_{01}$, $\phi_{02}=\phi_{01} + \phi_{12}$.
  Then
  \begin{gather}
    \label{eq:triciclo}
    \Fo_{\om_0}(\phi_{02}) = \Fo_{\om_0}(\phi_{01})
    +\Fo_{\om_1}(\phi_{12}).
  \end{gather}
\end{lemma}
(Same proof as in \cite[pp. 60f]{tian-libro}.)
\begin{lemma}[\protect{\cite[p. 59]{tian-libro}}]
  If $\phi_t $ is a differentiable family of smooth functions on $\XD$
  then
  \begin{align}
    &\desudt J_\om(\phi_t) = \frac{1}{\langle [\om]^n, [X]\rangle}\int_X
 \dphi_t
    \bigl ( \om^n -
    \om_t ^n \bigr)\\
    \label{eq:derivata-di-F0-lungo-la-curva}
    &\desudt \Fo_\om (\phi_t) = -\frac{1}{\langle [\om]^n, [X]\rangle}\int_X
    \dphi_t
    \om_t^n
  \end{align}
\end{lemma}

Assume now that $\om$ is a \Keler orbifold metric in the canonical
class, that is $\om\in 2\pi\chern_1(\xd)$.  Let $f=f(\om)\in
\cinf(\xd)$ be the unique function such that
\begin{equation}
  \label{eq:def-di-f}
  \Ric(\omega) - \omega = \idd f(\om)
  \qquad \int_X e^{ f(\om)}  = \int_X \om^n.
\end{equation}
Put $V=\MV=n!\vol(X)$ and define $A_\om, F_\om : \cinf\XD \ra \R$ by
\begin{gather}
  A_\om(\phi) = \log \biggl [ \intv e^{f(\om) - \phi}\om^n \biggr]
  \qquad F_\om(\phi ) = \Fo_\om(\phi) - A_\om(\phi).
\end{gather}
Using the notation of Lemma \ref{lemma-triciclico} if $\om_0, \om_1$
and $\om_2$ are \Keler metrics, then
\begin{gather}
  \label{eq:triciclo-per-F}
  F_{\om_0}(\phi_{02}) = F_{\om_0}(\phi_{01}) +F_{\om_1}(\phi_{12}).
\end{gather}
For $G\subset \Aut(\xd)$ a subgroup of isometries of
$(\xd,\om)$ put
\begin{equation}
  P_G(\xd, \om) =\{ \phi \in \cinf\XD : \om_\phi >0, \text{ and $\phi$ is
    $G$-invariant}\}. 
\end{equation}
If $G=\{1\}$ we simply write $P(\xd,\om)$.

In order to construct a \KE metric on $\XD$ the continuity method is
applied: fix a \Keler metric $\om $ in the canonical class and
consider the well-known equations
$$
(\om + \idd \phi_t)^n = e^{f - t\phi_t} \om^n \eqno{(*)_t}
$$
for a smooth family of functions in $\cinf\XD$.  Yau's estimates hold
for orbifold metrics, and in particular the Calabi conjecture is true,
which implies that \mc[0] admits a unique solution.  Denote by $\Lap$
the negative definite $\debar$-Laplacian on functions (that is  $ \Lap = -
\debar^* \debar$)  and by $-\lamma_j$ its eigenvalues.
\begin{lemma} [\protect{\cite[Theorem 4.20 p. 116]{aubin-some-problems}}]
  Let $\om$ be a \Keler metric on the compact orbifold $\XD$. If
  $\Ric(\om)\geq \eps >0$, then $\lamma_1\geq 1$.
\end{lemma}
It follows that the times $t$ for which \mc is solvable form an open
subset $S\subset [0,1]$ and that solutions $\phi_t$ are smooth in $t$,
see \cite[pp.  63-66]{tian-libro}. Given a $C^0$-estimate for the
solutions, Yau's estimates ensure that $S$ is closed, thus yielding
the solution up to $t=1$, which is a \KE metric.

\begin{prop}\label{stima-1}
  Let $\phi_t$ be a solution to \mc for $t\in [0, T_0)$.  Then
  $I_\om(\phi_t) - J_\om(\phi_t)$ is nondecreasing and
  $F_\om^0(\phi_t)\leq 0$.
\end{prop}
\dimo Differentiating \mc with respect to $t$ one gets
\begin{equation}
  \label{eq:derivata-dphi}
  (\Delta_t + t) \dot{\phi}_t = -
  \phi_t.  
\end{equation}
Therefore
\begin{multline*}
  \desudt \Bigl ( I_\om(\phi_t) - J_\om(\phi_t) \Bigr)
  =  \\
  =\intv \phi _t ( \phi_t + t \dot{\phi}_t) \om_t^n = (1-t^2) \intv
  \phi_t^2 \om_t ^n + \frac{1}{V}\int_X | \debar \dot{\phi}_t|^2
  \om_t.
\end{multline*}
This gives the first result. For the second use
\eqref{eq:derivata-di-F0-lungo-la-curva} and \eqref{eq:derivata-dphi}:
\begin{equation*}
  \begin{gathered}
    \desudt t \Fo_\om(\phi_t) = \Fo_\om(\phi_t)
    -\frac{t}{V}\int_X \dphi_t \om_t^n =  \\
    = \Fo_\om(\phi_t) +\frac{1}{V}\int_X (\Lap_t \dphi_t + \phi_t)
    \om_t^n = J_\om(\phi_t).
  \end{gathered}
\end{equation*}
Since $J_\om \geq 0$ the result follows.  \fine

The following estimates depend on the uniform Sobolev embedding (Lemma
\ref{uniform-Sobolev}) and their proof uses Moser iteration.
\begin{teo}
  [\protect{\cite[p. 67ff]{tian-libro}}] If $\phi_t$ is a family of
  solutions to \mc on the time interval $[0,T_0)$, then there is a
  constant $C=C(T_0)>0$ such that for any $t<T_0$
  \begin{gather}
    ||\phi_t||_\infty \leq C\bigr(1 + J_\om(\phi_t)\bigr )
    \label{eq:stimo-c0}
    \\
    0\leq -\inf_X \phi_t \leq C \biggl ( \intv (-\phi_t) \om_t^n +C
    \biggr )
    \\
    F_\om(\phi_t) \leq - A_\om(\phi_t) \leq C (1-t) \leq C.
    \label{eq:stimo-F-continuita}
  \end{gather}
\end{teo}

\begin{lemma} [\protect{\cite[\S 6]{bando-mabuchi-uniqueness}}]
  \label{good-gauge}
  Let $\XD$ be a Fano orbifold, $\omK$ a \KE metric and $\om$ a metric
  in the canonical class.  Then there is $g\in \Aut(\xd)$ such that $
  \om=g^*\omK + \idd \psi $ with $ \psi$ orthogonal to $\ker (
  \Delta_{g^*\omK} + 1)$ in $L^2(X,\om^n_{KE})$.
\end{lemma}

\begin{prop} [\protect{\cite[Prop. 5.3]{tian-KE97}}]
  \label{backward}
  Let $\XD$ be a Fano orbifold and $\omK$ a \KE metric in the
  canonical class. If $\om=\omK + \idd \psi$ is a \Keler metric, with
  $\psi \perp \ker(\Delta_{KE} +1 )$ and $\int_X e^{-\psi} \omK^n=0$,
  there is a solution $\{\phi_t\}_{t\in [0,1]}$ of \mc[t] with
  $\phi_0=0$ and $\phi_1=-\psi$.
\end{prop}

\begin{teo} [Ding-Tian]
  \label{ding-tian}
  If a Fano orbifold $\XD$ admits a \KE metric $\omK$, then $F_\om$ is
  bounded from below on $P(\X,\om)$ for any $\om$ in the canonical
  class.
\end{teo}
\dimo Thanks to \eqref{eq:triciclo-per-F} it is enough to bound
$F_\omK$. Given $\phi\in P(\xd, \omK)$ put $\om=\omK + \idd \phi$ and
let $g$ and $\psi$ be as in Lemma \ref{good-gauge}.  Using again
\eqref{eq:triciclo-per-F} it is enough to bound $F_{g^*\omK}(\psi)$.
Take a path as in Lemma \ref{backward}.  Thanks to Proposition
\ref{stima-1} $F_{g^*\omK} (\psi) = - F_\om(-\psi) = -F_\om (\phi_1) =
\Fo_\om(\phi_1) \geq0$.  \fine
\begin{remark}
  These estimates are enough to prove one half of Tian's fundamental
  theorem, namely that properness of $F_\om$ implies the existence of
  a \KE metric (see \cite[p. 63]{tian-libro}).
\end{remark}
The following normalisation of potentials is useful:
\begin{equation}
  Q_G(\X, \om) = \{ \phi \in P_G(\X, \om) : A_\om(\phi) = 0\}. 
\end{equation}
For any $\phi\in \PG$, $\phi +A_\om(\phi) \in Q_G(\X,\om)$.
\begin{prop}
  \label{criterio-del-sup}
  Let $\XD$ be a Fano orbifold, $\om\in 2\pi \chern_1(\xd)$ a \Keler
  metric and $G$ a compact group of isometries of $(\X,\om)$.  If
  there are constants $C_1 , C_2 > 0$ such that
  \begin{equation}
    F_\om(\phi) \geq C_1 \sup_X\phi - C_2
    \label{eq:prop1}
  \end{equation}
  for any $\phi \in Q_G(\X, \om)$, then $\XD$ admits a \KE metric.
\end{prop}
\dimo Let $\phi_t$ be a solution of \mc on $[0,T_0)$.  Since $\phi_t +
A_\om(\phi_t) $$ \in$$ Q_G(\xd,\om)$
\begin{equation}
  \label{eq:criterio-primo-su-P}
  \begin{gathered}
    F_\om(\phi_t) = F_\om\bigl (\phi_t + A_\om(\phi_t)\bigr)
    \geq \\
    \geq C_1 \sup_X \bigl (\phi_t + A_\om(\phi_t) \bigr ) -C_2 =C_1
    \sup_X \phi_t + C_1 A_\om(\phi_t) -C_2
  \end{gathered}
\end{equation}
Using \eqref{eq:stimo-F-continuita}
\begin{gather*}
  C_1 \sup_X \phi_t \leq F_\om(\phi_t) - C_1 A_\om(\phi_t) + C_2 \leq
  C_3 + C_2 +C_1 C_3.
\end{gather*}
Hence $\sup_X \phi_t $ is uniformly bounded. But $\Fo(\phi_t)\leq 0$,
so $J_\om(\phi_t)\leq \Fo_\om(\phi_t) + \sup\phi_t $ is bounded and
\eqref{eq:stimo-c0} yields the required bound of the $C^0$ norm.
\fine
\begin{lemma}[\protect{\cite[Lemma 2.3]{arezzo-ghigi-pirola}}]
  \label{superlemma-alpha}
  Let $\XD$ be a Fano orbifold, and $\om \in 2\pi\chern_1(M)$ a \Keler
  metric.  Then for any $\beta >0$ there are constants $C_1, C_2 >0$
  such that for any $\phi\in Q(\xd,\om)$
  \begin{equation}
    \label{eq:log-limita-sup}
    \log \biggl [ \intv e^{- (1+\beta)
      \phi} \om^n \biggr ]  \geq C_1 \sup_X \phi - C_2. 
  \end{equation} 
\end{lemma}
\begin{cor}
  \label{criterio-del-beta}
  If there are constants $C_1, C_2>0$ and $\beta >0$ such that
  \begin{equation}
    F_\om(\phi) \geq C_1 \log \biggl [ \intv e^{ - (1+\beta)
      \phi} \om^n \biggr ] -C_2
  \end{equation} 
  for any $\phi \in Q_G(\X, \om)$, then $\XD$ admits a \KE metric.
\end{cor}

\section{Existence theorems}

A \enf{current} on an orbifold \XD is a collection of
$\galf[\phi]$-invariant currents on any uniformiser $(U,\phi)$,
satisfying the usual compatibility condition with respect to
injections of uniformisers. In case $X$ is smooth, orbifold
differential forms on \XD are \emph{more} than ordinary differential
forms on $X$.  By duality orbifold currents on $\XD$ are \emph{less}
than ordinary currents on $X$: they are the continuous functionals on
$\est^k(X)$ that can be extended to the larger space $\est^k(\xd)$.
For \emph{positive} $(p,p)$-currents there is no difference between
the two notions, since every positive current has measure
coefficients, and every orbifold differential form has continuous
coefficients.  If $\gamma$ is a continuous hermitian form on a compact
orbifold $\XD$, an orbifold \emph{\Keler current} is a closed positive
(orbifold) current $T$ of bidegree (1,1) such that for some positive
constant $c$, $T\geq c \gamma$ in the sense of orbifold currents, that
is $\langle T-c\gamma, \eta \rangle \geq 0 $ for any positive $\eta
\in \est^{n-1,n-1}(\xd)$. The definition does not depend on the choice
of $\gamma$, since $X$ is compact.

If $\XD$ is a Fano orbifold, $G \subset \Aut(\xd)$ is a compact
subgroup and $\om$ is a $G$-invariant \Keler form in $2\pi
\chern_1(\xd)$, put
\begin{equation*}
  P^0_G(\xd,\om) = \{ \chi \in C^0(X) : \om + \idd \chi \text { is a
    \Keler orbifold current}\}. 
\end{equation*}
\begin{prop}
  \label{prop-estensione}
  (a) Any $\chi \in \barPG$ is the $C^0$-limit of a sequence
  $\phi_n\in P_G(\xd,\om)$. (b) The functionals $I_\om$, $J_\om$,
  $\Fo_\om$ and $F_\om$ can be extended to $\barPG$ and the extensions
  are continuous with respect to the $C^0$-topology.
\end{prop}
(See Prop. 2.2 and 2.3 in \cite{arezzo-ghigi-pirola}.)
\begin{lemma} [\protect{\cite[Lemma 2.6]{arezzo-ghigi-pirola} }]
  \label{Keler-current-image}
  If $\pi : (\xd_X) \ra \Y$ is an orbifold map between compact
  orbifolds, the direct image $\pi_*T$ of a \Keler current $T$ on
  $\XD$ is a \Keler current on $\Y$.
\end{lemma}
\dimo First of all observe that if $f : (X,\Delta_X) \ra (Y,\Delta_Y)$
is an orbifold map of degree $d$ and $\alf \in \est^{2n}(Y,\Delta_Y)$,
then $ \int_X f^*\alf = d\cdot \int_{Y} \alf.  $ Next let $\gamma_X$
and $\gamma_Y$ be continuous hermitian forms on $(\xd_X)$ and
$(Y,\Delta_Y)$ respectively.  Since $\pi^*\gamma_Y$ is continuous and
$\gamma_X$ is positive definite, there is $c_1 > 0$ such that
$\gamma_X \geq c_1 \pi^* \gamma_Y$.  If $T$ is a \Keler current on
$\XD$, by definition $T\geq c_2 \gamma_X$ for some $c_2 >0$, so $T\geq
c \pi^*\gamma_Y$ with $c=c_1 c_2>0$.  We want to prove that for any
positive form $\eta \in \est^{n-1,n-1}\Y$, $ \langle \pi_*T, \eta
\rangle \geq c\cdot \deg\pi \cdot \langle \gamma_Y, \eta \rangle$.
Choose orbifold charts $(V,\psi)$ on $\Y$ and $(U_i,\phi_i)$ on $\XD$
such that $\pi^\menuno\bigl( \psi(V)\bigr ) = \sqcup_i \phi_i(U_i)$.
Denote by $\tilde{T}_i$, $\tilde{\eta} $ and $\tilde{\gamma}_Y$ the
local representations in the orbifold charts and by $\tilde{\pi}_i :
U_i \ra V$ the liftings of $\pi$.  We can assume $\supp(\eta)\subset
\psi(V)$.  Then
\begin{equation*}
  \begin{gathered}
    \langle \pi_*T, \eta \rangle = \langle T, \pi^*\eta \rangle =
    \sum_i \frac{ \langle \tilde{T}_i, \tilde{\pi}^*_i
      \tilde{\eta}\rangle }{|\galf[\phi_i]|} \geq \sum_i \frac{c \cdot
      \langle \tilde{\pi}^*_i \tilde{\gamma}^*_Y, \tilde{\pi}^*_i
      \tilde{\eta}\rangle }{|\galf[\phi_i]|}
    =\\
    = \sum_i \frac{c}{|\galf[\phi_i]|} \int_{U_i} \tilde{\pi}^*_i
    \bigl (\tilde{\gamma}_Y \wedge \tilde{\eta} \bigr) = c \cdot\biggl
    ( \sum_i \frac{ \deg \tilde{\pi}_i}{|\galf[\phi_i]|} \biggr )
    \cdot \int_{V} \bigl (\tilde{\gamma}_Y \wedge \tilde{\eta} \bigr).
  \end{gathered}
\end{equation*}
Since
\begin{equation*}
  \sum_i \frac{  \deg \tilde{\pi}_i}{|\galf[\phi_i]|} = \frac{\deg
    \pi }{|\galf[\psi]|} 
\end{equation*}
we finally get
\begin{equation*}
  \langle \pi_*T, \eta \rangle  \geq c
  \int_{\psi(V)}  \bigl (\gamma_Y \wedge\eta \bigr) 
\end{equation*}
and this proves the lemma.  \fine{}
\begin{lemma}[\protect{\cite[Lemma 2.7]{arezzo-ghigi-pirola}}]
  \label{F0-nei-rivestimenti}
  Let $\pi: \XD \ra \Y$ be an orbifold map between $n$-dimensional
  \Keler orbifolds. Let $\om_Y$ be a \Keler metric on $\Y$ and $\chi
  \in P^0(Y,\Delta_Y,\om_Y)$ a \emph{continuous} potential such that
  $\pi^*\chi\in \cinf \XD$. Then
  \begin{equation}
    F^0_{\pi^*\om_Y} (\pi^*\chi)= F^0_{\om_Y}(\chi). 
  \end{equation}
\end{lemma}
\begin{teo}\label{one-covering}
  Let $(\xd_X)$ and $\Y$ be Fano orbifolds, $\pi : \XD \ra \Y$ an
  orbifold Galois covering of degree $d$ with $G= \galp$, $\om_Y$ a
  \KE metric on $\Y$ and $\om \in 2\pi\chern_1(\xd)$ a $G$-invariant
  \Keler metric.  Assume that numerically $ R^\orb(\pi) \equiv -\beta
  ( K_X + \Delta_X) $ for some $\beta \in \mathbb{Q}_+$.  Then there
  is a constant $C$ such that for any $\phi \in P_G(\xd, \om)$
  \begin{equation}
    F^0_{\om} (\phi) \geq  \frac{1}{1+\beta}
    \log \biggl [ \intv e^{- (1+\beta)\phi} \pi^*
    \om^n_Y \biggr ] -C. 
  \end{equation}
\end{teo}
The proof is identical to that of Theorem 2.2 in
\cite{arezzo-ghigi-pirola} and depends on the previous lemmata. Notice
that a $G$-invariant orbifold \Keler metric $\om$ always exists since,
according to Definition \ref{orbifold-Galois}, $G\subset \Aut\XD$.

\begin{teo}\label{criterio-del-c}
  Let $\XD$, $\XDi[1]$, $...,$ $\XDi[k]$ be $n$-dimensional Fano
  orbifolds. Assume that each $\XDi$ admits a \KE metric and that
  $\pi_i : \XD \ra \XDi$ are orbifold Galois coverings such that
  \begin{enumerate}
  \item the groups $\galf[\pi_i]$ are all contained in some
    \emph{compact} subgroup of $ \Aut\XD$;
  \item $ R^\orb(\pi_i) \equiv -\beta_i (K_X+ \Delta)$ for some
    $\beta_i\in \mathbb{Q}_+$.
  \end{enumerate}
  Define $\eta\in \cinf\XD$ by
  \begin{equation}
    \label{eq:def-di-eta}
    \frac{1}{k}\sum_{i=1}^k \pi_i^*\om_i^n = \eta \, \om^n,
  \end{equation}
  put $ c\perdef \sup\{\lambda \geq 0 : \eta^{-\lamma} \in L^1(X,
  \om^n)\} $ and $\beta \perdef \min \beta_i$. If
  \begin{equation}
    \label{eq:ipotesi-su-beta-c}
    \frac{1}{c} < \beta,
  \end{equation}
  then $\XD$ admits a \KE metric.
\end{teo}
The proof is the same as that of Theorem 2.3 and Proposition 2.4 in
\cite{arezzo-ghigi-pirola}.
\begin{remark}
  If there is only one covering ($k=1$) and $X$ is smooth, then $c$ is
  simply the \enf{complex singularity exponent} (that is the \enf{ log
    canonical threshold}) of the pair $(X,R^\orb)$ (see
  \cite{demailly-kollar-exponent} and
  \cite{kollar-singularities-pairs}).  On the other hand if there are
  enough coverings and the intersection of the ramification divisors
  $R^\orb(\pi_i)$ is empty, then $c=+\infty$ and
  \eqref{eq:ipotesi-su-beta-c} is automatically satisfied.
\end{remark}

\section{Applications}

Here we exhibit some concrete examples where
Theorem \ref{criterio-del-c} can be used to
prove the existence of \KE metrics on orbifolds.

\begin{teo}
  Let $X$ be a Fano manifold, $\sum_{i=1}^N D_i$ a divisor with local
  normal crossing and $\om$ a \KE metric on $X$. Given 
  integers $m_i>1$ put $\Delta = \sum_i (1-1/m_i)D_i$.  If $ \Delta
  \equiv - \delta K_X$ with $\delta \in (0,1)$ and
  \begin{equation}
    \label{eq:condizione-m-delta}
    m_i -1 < \frac{\delta}{1-\delta}
  \end{equation}
  for any $i=1, ..., N$, then $\XD$ is a Fano orbifold and has an
  orbifold \KE metric.
\end{teo}
\dimo $\XD$ is a Fano orbifold because $K_X + \Delta = (1-\delta )
K_X$ and $\delta<1$.  As observed in Example \ref{snc-galois} the map
$\id : \XD \ra X$ is an orbifold Galois cover and we want to apply
Proposition \ref{criterio-del-c} to it. The ramification divisor is
just $R^\orb=\Delta$ so
\begin{equation*}
  R^\orb(\id) = -\beta (K_X + \Delta)
\end{equation*}
with $ \beta = \delta/(1-\delta).  $ It remains to check that
\eqref{eq:condizione-m-delta} implies \eqref{eq:ipotesi-su-beta-c}.
Let $x$ be any point in $X$. Choose a system of coordinates $(V, z^1,
..., z^n)$ on $X$ as in Example \ref{local-normal-crossing} and let
$(U,\phi)$ be the corresponding orbifold chart for $\XD$ as in
\eqref{eq:snc-chart}.  Then on $\phi(U)=V$
\begin{equation}
  R^\orb=\Delta = \sum_{j=1}^k \Bigl( 1 - \frac{1}{m_j'} \Bigr ) \{z_j=0\}
\end{equation}
so that in the notation of \eqref{eq:def-di-eta}, $\eta (z) =
\gamma(z) |f(z)|^2$ on $U$, where $ f(z)= z_1^{m_1-1} \cdot \cdot
\cdot z_k^{m_k-1} $ and $\gamma$ is a smooth positive function.  Set $
c_x = \sup \{\lamma \leq 0: \int_U |f|^{-2\lamma} < +\infty\}.  $
Since
\begin{equation}
  \int_{U}  |f|^{-2\lamma}= \mathrm{const} \cdot
  \prod_{j=1}^k \int_{D} |z|^{-2\lamma (m'_j -1)} 
\end{equation}
where $D$ is the disk in $\C$, we get that $|f|^{-2\lamma} \in
L^1_{loc} $ on $U$ iff $ \lamma < 1/(m'_j -1).  $ So $ c_x =
\min\{1/(m'_j -1 ) : 1\leq j \leq k \}$,
\begin{equation}
  c=\sup_{x\in X} c_x = \min_i
  \frac{1}{m_i -1} 
\end{equation}
and
\begin{equation}
  \frac{1}{c} = \max (m_i - 1) < \frac{\delta}{1-\delta} =\beta.
\end{equation}
\fine
\begin{example}
  Let some divisors $D_i\in | \OO_{\PP^n}(d_i)|$, and some 
  integers $m_i>1$ be given for $i=1,...,N$. Let $m_1$ be the greatest
  of the $m_i$'s.  Put $\Delta = \sum_i (1-1/m_i)D_i$ and
  \begin{equation}
    \delta = \frac{\sum_i d_i \bigl (1 - \frac{1}{m_i} \bigr )} {n+1}. 
  \end{equation}
  Assume that
  \begin{enumerate}
  \item $\sum_i D_i$ is local normal crossing;
  \item $\delta < 1$;
  \item $m_1 (1-\delta) < 1$.
  \end{enumerate}
  Then the orbifold $\XD=(X,\Delta)$ admits an orbifold \KE metric of
  positive scalar curvature.
\end{example}

\begin{example}[Compare \protect{ \cite[Note
    36]{boyer-galicki-kollar-Annals}}]\label{new.exmps}
  Let $D_i$ be $n+2$ hyperplanes in general position in $\PP^n$:
$ D_i=\{z_i=0\} $ for $
  i=0,...,n, $ $ D_{n+1}=\{ z_0+...+z_n=0\}$. Set
$$
 \Delta= \sum_{i=0}^{n+1} (1-\tfrac{1}{m_i}) D_i.
$$
Then $(\PP^n, \Delta)$ has an orbifold \KE metric as soon as
\begin{equation}
  1< \sum_{i=0}^{n+1}\frac{1}{m_i} < 1 + (n+1) \min_i
  \frac{1}{m_i} 
  \end{equation}

  As in \cite{boyer-galicki-kollar-Annals}, many numerical examples
  come from Euclid's or Sylvester's sequence (cf.\
  \cite[A000058]{sequences}).  This is defined by the recursion
  relation
\begin{equation*}\label{extremseq}
c_{k+1}=c_1\cdots c_k+1=
c_k^2-c_k+1
\end{equation*}
beginning with $c_1=2$.
The sequence grows doubly exponentially, and it starts as
$$
2,3,7,43,1807, 3263443, 10650056950807,...
$$
It is easy to see  that 
$$
\sum_{i=1}^n\frac1{c_i}=1-\frac1{c_{n+1}-1}= 1-\frac{1}{c_1\cdots c_n}.
$$
We  get many new examples by taking
$$
(m_0=c_1, m_1=c_2, \dots,m_n=c_{n+1}-2,m_{n+1}).
$$
Then
$$
\sum_{i=0}^n \frac1{m_i}=1+\frac1{(c_{n+1}-1)(c_{n+2}-2)}.
$$
Thus our conditions are satisfied as long as
$$
c_{n+1}-2<m_{n+1}<n(c_{n+1}-1)(c_{n+2}-2)
$$
and $m_{n+1}$ is relatively prime to the other $m_i$.
\end{example}

\medskip

Another case when Theorem \ref{criterio-del-c} works is for degree 2
Del Pezzo surfaces $S$.  Here we consider the case when $S$ is allowed
to have cyclic quotient singularities. These are necessarily of the
form $\C^2/\Zeta_n$ where the group action is given by $(u,v)\mapsto
(\epsilon u,\epsilon^{-1}v)$ where $\epsilon$ is a primitive $n$-th
root of unity.  The $\Zeta_n$-invariant fuctions are generated by
$u^n,v^n,uv$.  This singularity is denoted by $A_{n-1}$.

For any degree 2 Del Pezzo surface $S$ the anticanonical class is
ample and it gives a degree 2 cover $\pi:S\to \PP^2$.  If $H$ denotes
the hyperplane class on $\PP^2$, then $-K_S=\pi^*H$.  The double cover
$\pi$ ramifies along a quartic curve $C$, thus
$R=\frac12\pi^*C=\pi^*2H$, $\beta=2$ and to apply Theorem
\ref{criterio-del-c} we need to ensure that $\eta^{-\lambda}$ be
integrable for $\lambda\leq \frac12$.  The singularities of $\pi$ lie
over the singularities of $C$, an $A_{n-1}$--singularity of $S$ lies
over an $A_{n-1}$--singularity of $C$ 
(cf.\ \cite[p.87]{barth-peters-vandeven}) 
and we can find local coordinates $(x,y) $
on $\PP^2$ such that $S$ is locally isomorphic to some neighbourhood
of the origin to the affine surface $ \{(x,y,t)\in \C^3: t^2=x^2 + 4
y^n\}$, the map $\pi$ being given simply by $\pi(x,y,t)=(x,y)$. An
orbifold chart is given by $\phi : U\subset \C^2 \ra S$ where
$\phi(u,v) = (u^n -v^n , uv, u^n + v^n)$. Thus $ \phi^* \pi^*
(dx\wedge dy)=n(u^n+v^n)\cdot du\wedge dv$ and
$\eta(u,v)=\textrm{const}\cdot|u^n+v^n|^2$.  It is easy to see by
direct integration or by blowing up (see e.g. \cite[Prop. 6.39 p.
168]{kollar-smith-corti}) that for $n\geq 2$, $|u^n+v^n|^{-2\lambda}$
is integrable if and only if $\lambda<\frac2{n}$. Thus Theorem
\ref{criterio-del-c} applies as long as $\frac12<c=\frac2{n}$, that is
for $n<4$. This proves Theorem \ref{teo-dp}.

One can also give a  different  proof of
the following result of  Mabuchi and Mukai \cite
[Corollary C]{mabuchi-mukai}.
\begin{teo}
  A diagonalizable singular  Del Pezzo surface of degree 4 
   admits an orbifold \KE metric.   
\end{teo}
A quartic Del Pezzo surface $S$ is the intersection of two quadrics in $\PP^4$,
$S=Q_1 \cap Q_2$. It is said to be diagonalizable if both $Q_1$ and
$Q_2$ can be put simultaneously in diagonal form. 
If $S$ is singular then in suitable coordinates
it is given by equations
\begin{equation*}
    h_0:= 
 x_0^2 + x_1 ^2 + x_2^2
    +x_3^2 +x_4^2 =0\quad \mbox{and}\quad 
    h_1:=
  \lamma_2 x_2^2
    +\lamma _3 x_3^2 +\lamma_4 x_4^2 =0
\end{equation*}
If two of the  $\lamma_i$ coincide 
then $S$ is a quotient of
$\PP^1 \times \PP^1$ and so has an orbifold \KE metric (see
\cite[p.136]{mabuchi-mukai}). Thus assume that the
 $\lamma_i$ are distinct nonzero complex numbers. 
 For $i=2,3,4$, 
the equation $\lamma_i h_0-h_1=0$ does not involve $x_i$,
and by dropping the $x_i$ variable we get 
 smooth quadrics
$$
Q_i=\{ (\lamma_i h_0-h_i=0)\}\subset \PP^3.
$$
The map $\pi_i : S \ra Q_i$ given by forgetting $x_i$ is a double cover
ramified over the hyperplane section $S\cap\{x_i=0\}$. Since the $Q_i$
are smooth two-dimensional quadrics, they are
K\"ahler-Einstein. 
On the other hand, the divisors $R^\orb (\pi_i) $ are disjoint, 
so $\eta$ is strictly positive on all $S$,   $c=\infty$ and  Theorem
\ref{criterio-del-c} yields that $S$ admits an orbifold \KE metric.

\begin{ack}  
We thank C.\ Arezzo, Ch.\ Boyer, K.\ Galicki and G.P.\ Pirola for
useful comments.
J.K.\ 
was partially supported by the NSF under grant number
DMS-0200883. 
\end{ack}

\vskip.3cm

\noindent Universit\`a di Milano Bicocca
\begin{verbatim}alessandro.ghigi@unimib.it\end{verbatim}

\noindent Princeton University, Princeton NJ 08544-1000
\begin{verbatim}kollar@math.princeton.edu\end{verbatim}

\end{document}